\documentclass{article}
\usepackage{amsmath, amsthm, amssymb, epsfig}
\author{Claire Levaillant}
\title{On arithmetical structures on $K_9$}

\newcommand{\nts}{\negthickspace}
\newcommand{\lra}{\longrightarrow}
\newcommand{\ms}{\mathcal{S}}
\newcommand{\h}{\widehat{2}}
\newcommand{\mbf}{\mathbf}
\newcommand{\xk}{x^{'}_k}
\newcommand{\xl}{x^{'}_l}
\newcommand{\xu}{x^{'}_1}
\newcommand{\xd}{x^{'}_2}
\newcommand{\xt}{x^{'}_3}
\newcommand{\xq}{x^{'}_4}

\newcommand{\xit}{x^{'}_{i_t}}
\newcommand{\xis}{x^{'}_{i_s}}
\newcommand{\xiu}{x^{'}_{i_1}}
\newcommand{\xiku}{x^{'}_{i_{k_1}}}
\newcommand{\xiksmu}{x^{'}_{i_{k_{s-1}}}}
\newcommand{\xiv}{x^{'}_{i_v}}
\newcommand{\al}{\alpha}
\newcommand{\be}{\beta}
\newcommand{\ga}{\gamma}
\newcommand{\lb}{\lbrace}
\newcommand{\rb}{\rbrace}

\begin{document}
\maketitle
\begin{abstract}
\noindent We study the arithmetical structures on the complete graph $K_9$. \\Our method is based on studying the solutions to writing the unit as a sum of 9 unit fractions. We work from the perspective of the diophantine equation and use some elementary properties on the $p$-adic valuations. The proofs are assisted by trees and automata.
\end{abstract}

\section{Historical background}
Egyptian fractions appear in an ancient Egyptian text, namely the Rhind Mathematical Papyrus, dating back to around $1650$ B.C. \cite{ROB}. The fact that any positive rational number $r$ has an Egyptian fraction representation of the form
$$r=\frac{1}{x_1}+\dots+\frac{1}{x_n}\qquad\text{with $1\leq x_1<\dots <x_n$}$$
is due to Fibonacci in $1202$, using a greedy algorithm (at each step the algorithm chooses greedily the largest possible unit fraction that can be used in any representation of the remaining fraction) and later rediscovered by Sylvester in \cite{SYL}. Fibonacci's greedy algorithm consists of repeatedly performing:

$$\frac{x}{y}=\frac{1}{\lceil \frac{y}{x}\rceil}+\frac{(-y)\,\text{mod}\,x}{y\lceil \frac{y}{x}\rceil}$$

Later on, several authors tried to improve the Fibonacci-Sylvester algorithm, for instance by minimizing the number of unit fractions.
A non exhaustive list is \cite{ERD2}\cite{GOL}\cite{BLE}\cite{BER}. Recently, Louwsma and Martino studied the rational numbers with odd greedy expansion of fixed length. Their work appears in \cite{LOM}.

Problems and results concerning Egyptian fractions appear in \cite{ERD1} and \cite{BLE}. \\
The present paper is concerned with representations of the unit, that is setting $r=1$ and the $x_i$'s are not necessarily distinct. We present a new viewpoint for the arithmetic study of the diophantine equation, based on elementary properties on $p$-adic numbers. One interest in the solutions to this diophantine equation relies in the fact that finding such solutions to the diophantine equation is equivalent to finding $n$ positive integers $k_i$'s with $1\leq i\leq n$ with no common factor such that
$$k_j\bigg|\sum_{i=1}^n k_i\qquad\text{for all $j$}$$
This problem is itself equivalent to finding an arithmetical structure on the complete graph $K_n$, see for instance Theorem $4.4$ of \cite{KER} proven in a more general setting. \\
Arithmetical structures were originally introduced by Lorenzini in \cite{LOR} in terms of matrices as a generalization of the Laplacian matrix in order to study intersections of degenerating curves in algebraic geometry. If $G$ is a finite connected graph with $n$ vertices and $A$ is the adjacency matrix of $G$, then an arithmetical structure on $G$ is a pair of vectors $(\mathbf{d},\mathbf{r})$ whose entries are nonnegative and positive integers respectively and the entries of $\mathbf{r}$ have no nontrivial common factor, such that:
$$\big(diag(\mathbf{d})-A\big)\mathbf{r}=\mathbf{0}$$
The critical group associated with an arithmetical structure \cite{KAP} generalizes the sandpile group in the special case of the Laplacian matrix which is by definition the difference of the diagonal degree matrix with the adjacency matrix. The sandpile group itself relates to the chip-firing game, see e.g. \cite{GLA} and \cite{BLS}.
Lorenzini showed in his Lemma $1.6$ of \cite{LOR} that any finite, connected graph has a finite number of arithmetical structures. In \cite{MLA}, the authors show that the number of arithmetical structures on the path graph $P_n$ is given by the Catalan number $C(n-1)$ and that the number of arithmetical structures on the cycle graph $C_n$ is given by the binomial coefficient $\binom{2n-1}{n-1}$. In \cite{LEA}, the authors treat the case of graphs corresponding to Dynkin diagrams of type $D_n$. For a given $n\geq 4$, they establish the recursive formula:
$$|Arith(D_n)|=2C(n-2)+\sum_{m=4}^nB(n-3,n-m)|SArith(D_m)|,$$
where $C(n)=\frac{1}{n+1}\binom{2n}{n}$ is the Catalan number, $B(n,k)=\frac{n-k+1}{n+1}\binom{n+k}{n}$ is the ballot number and $|SArithm(D_m)|$ is the number of smooth arithmetical structures on $D_m$. The authors define a smooth arithmetical structure on $D_n$ by imposing that $d_i\geq 2$ at each vertex $i$, except possibly at the trivalent vertex say $i=t$ (in fact, they show in their Proposition $2.3$ that these conditions on the $d_i$'s imply $d_t=1$). They then come up with bounds for the number of smooth arithmetical structures and ultimately
deduce the following bounds for the arithmetical structures on $D_n$ with $n\geq 4$:
$$2C(n-2)+C(n-3)\leq |Arith(D_n)|<2\,C(n-2)+702\,C(n-3)$$
In \cite{GLW}, the authors deal with arithmetical structures on paths with a doubled edge and conjecture how the number of arithmetical structures grows depending on the path length and the location of the doubled edge. \\
The importance of arithmetical structures on complete graphs gets reinforced by Conjecture $6.10$ of \cite{COR} which asserts that for any connected simple graph $G$ with $n$ vertices, the number of arithmetical structures on $G$ is at most the number of arithmetical structures on $K_n$. The number of arithmetical structures on $K_n$ appears in the online encyclopedia of integer sequences \cite{OEI} for $n\leq 8$. In the present paper, we investigate arithmetical structures on $K_9$ from the perspective of the Diophantine equation, that is we look for solutions to $1/x_1+\dots+1/x_9=1$ in positive integers $x_i$'s, $i=1,\dots,9$. It is shown in \cite{BUR} that there are exactly five solutions in distinct and odd integers. The solutions in possibly non distinct and possibly even integers have not been investigated so far. \\
Our paper is structured as follows. We will do some general $p$-adic analysis on the Diophantine equation in $\S\,2$ and then apply it specifically to the case of $9$ positive integer variables in $\S\,3$. Like other authors also did it, we will restrict the study to some specific cases described in $\S\,3$. Within these specific cases, the method that we use allows to count the number of solutions for the general complete graph $K_n$ with $n\geq 9$. We wrote a Mathematica program which given the number $n$ in entry returns this number of solutions. Last, still in the framework of these specific cases, we wrote a program in CAML allowing to obtain the complete list of solutions and applied our program in the cases of $K_9-K_{13}$.\\

In the statements below and throughout the paper, $(E)_n$ denotes the Diophantine equation
\begin{equation}\frac{1}{x_1}+\dots+\frac{1}{x_n}=1\tag*{$(E)_n$}\end{equation}

Our main results are the following. Theorem $1$ deals with $n=9$ specifically. Theorem $2$ provides the count of the solutions for a general $n\geq 9$ using a Mathematica program. Theorem $3$ exhibits all the solutions when $n=10,11,12,13$ using a CAML program.

\newtheorem{Theorem}{Theorem}
\begin{Theorem}
Consider the Diophantine equation $(E)_9$ in distinct integers.
There are $54$ solutions to the equation of the form $2^{\al}q^{\beta}$ with $q$ odd prime and $\al\leq 2$. These are exhibited along the paper and listed as:
$X_1-X_{17}$, $Y_1$, $Y_2$, $Z_1-Z_{27}$, $\widehat{Z_1}$, $\widehat{Z_2}$, $\widehat{Z_3}$, $T_1$, $T_2$, $\widehat{T_1}$, $U$ and $V$.
\end{Theorem}

\begin{Theorem} Consider the Diophantine equation $(E)_n$ with $n\geq 9$ in distinct integers. The total number of solutions to $(E)_n$ of the form $2^{\al}q^{\be}$ with $q$ odd prime and $\al\leq 2$ gets provided by adding the first seven integers of the output of the Mathematica program of Appendix G plus two (resp plus one) in the case when $n$ is odd (resp $n$ is even).
\end{Theorem}

\begin{Theorem} The cases of $(E)_{10}$, $(E)_{11}$, $(E)_{12}$ and $(E)_{13}$ in distinct integers using a CAML program.\\
Consider the Diophantine equation $(E)_{n}$ with $n=10,11,12,13$ in distinct integers.
There are respectively $101$, $192$, $363$, $692$ solutions to the equation of the form $2^{\al}q^{\beta}$ with $q$ odd prime and $\al\leq 2$ whose respective $100$, $190$, $362$ and $690$ first solutions are all exhibited in Appendix E when $q=3$. The additional solutions are
$U_{10}-U_{13}$ when $q=5$ and $V_{11},\,V_{13}$ respectively corresponding to $n=11$ and $n=13$ when $q=7$, and these solutions get exhibited in $\S\,4$ of the paper.
\end{Theorem}

\section{Arithmetics of arithmetical structures on the complete graph $K_n$}
The goal of this part is to study for a given prime $p$ the $p$-adic valuations of the positive integers $x_i$'s when these constitute a solution to the equation $(E)_n$. We begin our study with the fundamental property stated right below.
\newtheorem{Proposition}{Proposition}

\begin{Proposition}
Let $p$ be a prime number. If $(x_1,\dots,x_n)$ is a solution, then $v_p(x_k)=v_p(x_l)$ for some $k\neq l$.
\end{Proposition}

\newtheorem{Lemma}{Lemma}

\begin{Lemma}
Let $p$ be a prime number. If $(x_1,\dots,x_n)$ is a solution, then $v_p(x_i)=0$ for some $i$ or $v_p(x_k)=v_p(x_l)$ for some $k\neq l$.
\end{Lemma}

\textsc{Proof of Lemma $1$}. By assumption, the $x_i$'s satisfy to
\begin{equation}\sum_{i=1}^r x_1\dots\hat{x_i}\dots x_n-x_1\dots x_n=0\end{equation}
The $p$-adic absolute value being ultrametric, for any $p$-adic numbers $a_1,\dots,a_s$, we have:
$$v_p(a_1+\dots+a_s)\geq Min_{1\leq i\leq s}v_p(a_i),$$
with the equality holding when the minimum is unique.\\
By convention, we have $v_p(0)=+\infty$. Then $(1)$ implies that:
\begin{equation}v_p(x_1\dots x_n)=v_p(x_1\dots\hat{x_i}\dots x_n),\qquad\text{some $1\leq i\leq n$}\end{equation}
or
\begin{equation}
v_p(x_1\dots\hat{x_i}\dots x_n)=v_p(x_1\dots\hat{x_j}\dots x_n),\;\;\;\,\text{some $1\leq i<j\leq n$}
\end{equation}
Moreover, if $(2)$ holds, then we get $v_p(x_i)=0$ and if $(3)$ holds, we get $v_p(x_i)=v_p(x_j)$.\hfill $\square$\\

\textsc{Proof of Proposition $1$.} If $p$ divides none of the $x_i$'s, then for all $i$ with $1\leq i\leq n$, we have $v_p(x_i)=0$ and so we are done.
From now on, suppose that $p$ divides at least one of the $x_i$'s.
The fact that $(x_1,\dots,x_n)$ is a solution forces $v_p(x_i)=0$ for some $i$ or $v_p(x_k)=v_p(x_l)$ for some $k\neq l$. In the latter case, we are done. Otherwise, assume without loss of generality that $v_p(x_1)=0$ and factor the diophantine equation as follows.
\begin{equation}
x_1(\sum_{i=2}^n x_2\dots\hat{x_i}\dots x_n-x_2\dots x_n)+x_2\dots x_n=0
\end{equation}
It comes:
\begin{equation}
v_p(x_2\dots x_n)=v_p(\sum_{i=2}^n x_2\dots\hat{x_i}\dots x_n-x_2\dots x_n)
\end{equation}
By assumption, there exists an integer $s$ amongst $\lbrace 2,\dots,n\rbrace$ such that $v_p(x_s)>0$. Then,
\begin{equation}
v_p(x_2\dots x_n)>Min_{2\leq i\leq r}v_p(x_2\dots \hat{x_i}\dots x_r)
\end{equation}
In turn, there must exist two integers $k$ and $l$ with $2\leq k<l\leq n$ such that:
\begin{equation}
v_p(x_2\dots\hat{x_k}\dots x_n)=v_p(x_2\dots\hat{x_l}\dots x_n)
\end{equation}
Then, for these two distinct integers $k$ and $l$, we have $v_p(x_k)=v_p(x_l)$ and moreover, these two valuations are non-zero. Hence the straightforward corollary.
\newtheorem{Corollary}{Corollary}
\begin{Corollary}
Let $p$ be a prime number. If $(x_1,\dots,x_n)$ is a solution and $p$ divides at least one of the $x_i$'s, then there exist two distinct integers $k$ and $l$ such that $v_p(x_k)=v_p(x_l)$ and this common valuation is non-zero.
\end{Corollary}
\begin{Corollary}
Let $p$ be a prime number. If $(x_1,\dots,x_n)$ is a solution, then the highest power of $p$ occurring as a factor amongst the $x_i$'s occurs at least twice.
\end{Corollary}
\textsc{Proof of Corollary $2$}. If $p$ divides none of the $x_i$'s, the result is trivial. Otherwise, by Corollary $1$, there exist distinct integers $k$ and $l$ with $v_p(x_k)=v_p(x_l)=\alpha>0$. Either $\al$ is the highest integer such that $p^{\al}$ occurs as a factor amongst the $x_i$'s, in which case we are done or we have:
\begin{equation}
Max_{j\not\in\lbrace k,l\rbrace}v_p(x_j)>\al
\end{equation}
Suppose for a contradiction that this maximum is unique and attained at $x_s$ with $s\not\in\lbrace k,l\rbrace$.\\
Write:
$$\begin{array}{cccc}x_k&=&p^{\al}\xk&\text{with $p\not|\xk$}\\
x_l&=&p^{\al}\xl&\text{with $p\not|\xl$}
\end{array}$$
After division by $p^{\al}$, the diophantine equation reads:
\begin{equation}
p^{\al}\xk\xl\prod_{j\not\in\lbrace k,l\rbrace}x_j=\xk\prod_{j\not\in\lbrace k,l\rbrace}x_j+\xl\prod_{j\not\in\lbrace k,l\rbrace}x_j+\sum_{m\not\in\lbrace k,l\rbrace}p^{\al}\xk\xl\prod_{j\not\in\lbrace m,k,l\rbrace} x_j
\end{equation}
It follows that:
\begin{equation}
v_p\big(\prod_{j\not\in\lb k,l\rb}x_j\big)\leq\al +\sum_{j\not\in\lb k,l,s\rb}v_p(x_j)
\end{equation}
This implies $v_p(x_s)\leq\al$, a contradiction. \hfill $\square$\\\\
In what follows, we denote the elementary symmetric functions of $a_1,\dots,a_n$ by $\sigma_i(a_1,\dots,a_n)$ with $1\leq i\leq n$. And so by definition,
$$\sigma_i(a_1,\dots,a_n):=\sum_{1\leq k_1<\dots<k_i\leq n}a_{k_1}\dots a_{k_i}$$
Using this notation, the following result holds.
\begin{Corollary}
Let $p$ be a prime number and $(x_1,\dots,x_n)$ be a solution, with $p$ dividing at least one of the $x_i$'s. Suppose the highest power of $p$, say $p^{\al}$, dividing the $x_i$'s occurs exactly at indices $i_1,\dots, i_s$ with $s\geq 2$ and set for $t\in\lb 1,\dots,s\rb$, $x_{i_t}=p^{\al}\xit$ with $p\not\!|\xit$. \\Then, $$Max_{j\not\in\lb i_1,\dots,i_s\rb}v_p(x_j)<\al$$ and
$$p^{\al-Max_{j\not\in\lb i_1,\dots,i_s\rb}v_p(x_j)}|\sigma_{s-1}(\xiu,\dots,\xis)$$
Moreover, if $Max_{j\not\in\lb i_1,\dots,i_s\rb}v_p(x_j)$ is unique and positive, then
$$v_p(\sigma_{s-1}(\xiu,\dots,\xis))=\al-Max_{j\not\in\lb i_1,\dots,i_s}v_p(x_j)$$
\end{Corollary}
\textsc{Proof of Corollary $3$}. We will write a similar equation as the one in $(9)$ with a larger number of $x^{'}_i$'s involved.
We have, after dividing the diophantine equation by $p^{(s-1)\al}$:
\begin{equation}\begin{split}
\sigma_{s-1}(\xiu,\dots,\xis)\prod_{j\not\in\lb i_1,\dots,i_s\rb}x_j&+p^{\al}\xiu\dots\xis\sum_{m\not\in\lb i_1,\dots,i_s\rb}\prod_{j\not\in\lb m,i_1,\dots,i_s\rb}x_j\\&=p^{\al}\xiu\dots\xis\prod_{j\not\in\lb i_1,\dots,i_s\rb}x_j
\end{split}\end{equation}
It follows immediately that:
\begin{equation}
v_p(\sigma_{s-1}(\xiu,\dots,\xis))\geq \al-Max_{m\not\in\lb i_1,\dots,i_s\rb}v_p(x_m),
\end{equation}
with equality holding when the maximum is unique and positive.\hfill $\square$\\

\begin{Corollary}
Let $p$ be a prime number and $(x_1,\dots,x_n)$ be a solution, with $p$ dividing at least one of the $x_i$'s. Suppose the highest power of $p$, say $p^{\al}$, dividing the $x_i$'s occurs exactly at indices $i_1,\dots, i_s$ with $s\geq 2$ and set for $t\in\lb 1,\dots,s\rb$, $x_{i_t}=p^{\al}\xit$ with $p\not\!|\xit$. Let $k_1,\dots,k_{s-1}$ be $(s-1)$ distinct integers chosen amongst $\lb 1,\dots,s\rb$.
Then,
$$v_p\big(\sigma_{s-2}(\xiku,\dots,\xiksmu)\big)=0$$
\end{Corollary}
\textsc{Proof of Corollary $4$}. Eq. $(11)$ still holds with a smaller number of $\xit$'s involved. In particular, by writing $(11)$ with only $(s-1)$ of the $s$ $\xit$'s, we see that the maximum that is involved is unique and equal to $\al$. It yields the result. \hfill $\square$\\
\begin{Corollary}
Let $p$ be a prime number and let $(x_1,\dots,x_n)$ be a solution with $p$ dividing at least one of the $x_i$'s. Suppose the highest power of $p$, say $p^{\al}$, occurs at indices $i_1,\dots,i_t$ with $t\geq 2$ and set for $v\in\lb 1,\dots, t\rb$, $x_{i_v}=p^{\al}\xiv$ with $p\not\!|\xiv$. Let $s$ be the highest such $t$.\\
If $p|\sigma_{t-1}(\xiu,\dots,\xit)$, then either $s=t$ or $s\geq t+2$.

\end{Corollary}
The latter three corollaries are best illustrated on an example.\\
When $n=11$, Burshtein proves that there are exactly $17$ solutions $\mathbf{B_1}-\mathbf{B_{17}}$ in distinct integers of the form $3^{\al}5^{\be}7^{\ga}$ and exhibits these solutions in \cite{BUR2}.
Consider for instance
$$\mathbf{B_{12}}=\lb 3,5,7,3^2,3^2.5.7,5^2,3^3,7^2,7^2.3^3,7.5^2.3^3\rb,$$
written in terms of prime numbers. Set $p=3$ and so $\al=3$ using the previous notations. Visibly $3^3$ occurs $3$ times as a factor.
Obviously, $3$ does not divide $\sigma_1(7.5^2,1)=176$, nor $\sigma_1(7^2,1)=50$, nor $\sigma_1(7^2,7.5^2)=224$, just like implied by Corollary $4$. The highest other $3$-adic valuation is $2$ and we have
$v_3(\sigma_2(7^2,7.5^2,1))=v_3(8799)\geq 1$, just like implied by Corollary $3$. Set now $p=5$ and so $\al=2$. Visibly, $5^2$ occurs exactly twice as a factor and indeed, we have $5|1+7.27\,(=190)$

\section{Combinatorics of arithmetical structures on the complete graph $K_9$ in some restricted cases.}

$\qquad$ The case $n=9$ studied here is special in that when the integers are odd, the minimal length of representation of the unit as an Egyptian fraction is $9$. This was first predicted by J. Leech, see \cite{GUY}, pp. $89$. In \cite{BUR}, Burshtein proved that the  diophantine equation has exactly $5$ solutions that were first identified by S. Yamashita in $1976$, see the online The prime puzzles and problems connection, Problem $35$ entitled "More wrong turns..." Ronald Graham had once asked Andr\'e Weil why the Egyptians did this, that is to represent fractions as a sum of unit fractions, and Weil's answer was "They took a wrong turn..." \\The five solutions in odd and distinct integers copied from \cite{BUR} are:
$$\left\lb\begin{array}{l}
\mathbf{B_1}=\lb 3,5,7,9,11,15,21,231,315\rb\\
\mathbf{B_2}=\lb 3,5,7,9,11,15,35,45,231\rb\\
\mathbf{B_3}=\lb 3,5,7,9,11,15,21,135,10395\rb\\
\mathbf{B_4}=\lb 3,5,7,9,11,15,33,45,385\rb\\
\mathbf{B_5}=\lb 3,5,7,9,11,15,21,165,693\rb\end{array}\right.$$

We are now concerned with searching for solutions when the $x_i$'s are not necessarily distinct and are possibly even.
First and foremost, we will search for solutions with $2$ occurring at least twice as a factor.

\begin{Lemma} Let $(x_1,\dots,x_9)$ be a solution and at least one of the $x_i$'s is even. The highest power of $2$ occurring as a factor amongst the $x_i$'s must occur an even number of times.
\end{Lemma}
\textsc{Proof.} Using our usual notations, set $\al$ the highest positive $2$-adic valuation amongst the $x_i$'s and let $s$ denote the number of occurrences. By Corollary $2$, we know that $2^{\al}$ occurs at least twice as a factor amongst the $x_i$'s. Without loss of generality, $x_1=2^{\al}\xu$ and $x_2=2^{\al}\xd$ with $\xu$ and $\xd$ both odd. The sum of two odd numbers being even, we have $2|\xu+\xd$, which by Corollary $5$ implies $s=2$ or $s\geq 4$. Suppose now $s\neq 2$. Without loss of generality,
$x_i=2^{\al}x^{'}_i$ for $i=1,\dots,4$ and each $x^{'}_i$ is odd. Now, $2$ divides $\sigma_3(\xu,\xd,\xt,\xq)$ since the latter elementary symmetric function is the sum of $4$ odd terms. A new application of Corollary $5$ yields $s=4$ or $s\geq 6$. The result follows from proceeding by induction, since there is an even number of ways of picking $2l-1$ elements amongst $2l$ elements and a sum of an even number of odd integers is an even integer.\hfill $\square$\\

For the remainder of our discussion, it will be convenient to denote by $\al_p$ the highest $p$-adic valuation of the $x_i$'s and by $s_p$ the number of occurrence(s) amongst the $x_i$'s. For instance, we have just shown that $s_2\in\lb 0,2,4,6,8\rb$. \\
Set also $X:=\lb x_1,\dots,x_9\rb$ as the set of constituents of a solution to the diophantine equation. If along the proofs, we need to refer to solutions with $k$ vertices with $k$ not necessarily equal to $9$, we will specify it by writing $X_{n=k}$ instead of $X$.\\\\
\textbf{In what follows, we assume that the $x_i$'s are distinct and $2$ is a divisor.}



\begin{Lemma}(i) Suppose $2\not\in X$ and $\al_2=1$. Then, there is no solution with only two primes.\\
\hspace{0.4cm}(ii) Suppose $2\in X$ and $\al_2=1$. Then,
$$\mathbf{Z}_1=\lb 2,3,3^2,3^3,3^4,3^5,3^6,3^7,2.3^7\rb$$
is the only solution with only two primes.\\
\end{Lemma}


\begin{Lemma}
$(i)$ If $2\not\in X$ and $\al_2=2$, then the solutions with only two primes are precisely $Y_1$, $Y_2$, $T_1$, $T_2$ and
\begin{eqnarray*}
\mathbf{\widehat{T}_1}&:=&\lb 2^2,3,2.3,3^2,2.3^2,2^2.3^2,3^3,3^4,2.3^4\rb
\end{eqnarray*}
$(ii)$ If $2\in X$ and $\al_2=2$, then the only solutions in $2^{\al}.q^{\beta}$ with $q\neq 3$ are:
$$\mathbf{V}:=\lb 2,2^2,7,2.7,7^2,2.7^2,7^3,2.7^3,2^2.7^3\rb$$
and
$$\mathbf{U}:=\lb 2,2^2,5,5^2,5^3,5^4,5^5,5^6,2^2.5^6\rb$$
The solutions in $2^{\al}.3^{\beta}$ are $\widehat{Z}_1-\widehat{Z}_3$, $Z_2-Z_{27}$ and $X_1-X_{17}$.
\end{Lemma}

We will first prove Lemma $4$ and then prove Lemma $3$ by using some elements in the proof of Lemma $4$. \\

\textsc{Proof of Lemma $4$.} In the proof, there is no distinction between $2\in X$ or $2\not\in X$. Having $2\not\in X$ is more rare, hence the distinction in the statement itself. We first lay out some general considerations based on $\al_2=2$ and our results of $\S\,2$. Let $q$ be an odd divisor. Our assumptions imply that the highest power of $q$ occurs twice or thrice. We apply Corollary $3$ to $q$.

If the highest power of $q$ occurs three times, then $q|\sigma_2(1,2,4)=14$. It forces $q=7$. Moreover, still by Corollary $3$, the second largest valuation must be $\al_7-1$. Further, we can show that this second largest valuation may only occur as $\lb 7^{\al_7-1},2.7^{\al_7-1}\rb$. Indeed, in all the other cases, by contracting $\lb 7^{\al_7},2.7^{\al_7},2^2.7^{\al_7}\rb$ to $2^2.7^{\al_7-1}$ and possibly contracting further, we would get a solution of $X_{n=k}$ for some adequate $k\leq 7$ with the highest $7$-valuation occurring less than three times. Continuing this process inductively, it forces $2\in X$ and moreover, we obtain the result of the first point of $(ii)$, with $s_2=2$ the only available option.

Suppose now the highest power of $q$ occurs exactly twice. Then $q|\sigma_1(1,2)=3$ or $q|\sigma_1(1,2^2)=5$ or $q|\sigma_1(2,2^2)=2.3$. This implies $q=3$ or $q=5$. \\In the first case, the last two terms of the sum are either $\big\lb \frac{1}{3^{\al_3}},\frac{1}{2.3^{\al_3}}\big\rb$ or $\big\lb\frac{1}{2.3^{\al_3}},\frac{1}{2^2.3^{\al_3}}\big\rb$ and the second largest $3$-valuation occurring is $\al_3-1$. \\In the second case, the last two terms of the sum are $\frac{1}{5^{\al_5}}$ and $\frac{1}{2^2.5^{\al_5}}$. Again, the second largest $5$-valuation occurring is $\al_5-1$.\\
We first deal with the case $q=5$. This time, it will be useful to note that $\lb 5^{\al_5},2^2.5^{\al_5}\rb$ contracts to $2^2.5^{\al_5-1}$. By the same reasoning as for $q=7$ above, the only option is to have $5^{\al_5-1}$ in the solution. We then work inductively until getting a solution of $X_{n=3}$ of the form $\lb x_1,x_2,2^2.5^{\al_5-6}\rb$ forcing $x_1=2$, $x_2=4$ and $\al_5=6$. Consequently the only possible original solution is $U$, provided in point $(ii)$ of the lemma.\\
When $q=3$, we treat both cases evoked above simultaneously and inductively on the number of vertices. \\

We state below a series of facts on which our discussion will be based.
\newtheorem{Fact}{Fact}
\begin{Fact} Suppose the second highest $3$-valuation $\al_3-1$ occurs exactly once.
Then,\\
(i) It is impossible to have $\lb 3^{\al_3-1},2^23^{\al_3},2.3^{\al_3}\rb\subseteq X$.\\
(ii) If $\lb 2.3^{\al_3-1},2^2.3^{\al_3},2.3^{\al_3}\rb\subseteq X$, then the solution of $X_{n=9}$ is a solution of $X_{n=8}$ ending in $\lb 2.3^{\al_3-1},2^2.3^{\al_3-1}\rb$.\\
(iii) It is impossible to have $\lb 2^2.3^{\al_3-1},2^2.3^{\al_3},2.3^{\al_3}\rb\subseteq X$.\\\\
(iv) If $\lb 3^{\al_3-1},3^{\al_3},2.3^{\al_3}\rb\subseteq X$, then the solution of $X_{n=9}$ is a solution of $X_{n=8}$ ending in $\lb 3^{\al_3-1},2.3^{\al_3-1}\rb$.\\
(v) It is impossible to have $\lb 2.3^{\al_3-1},3^{\al_3},2.3^{\al_3}\rb\subseteq X$.\\
(vi) If $\lb 2^2.3^{\al_3-1},3^{\al_3},2.3^{\al_3}\rb\subseteq X$, then the solution of $X_{n=9}$ is a solution of $X_{n=8}$ ending in $\lb 2.3^{\al_3-1},2^2.3^{\al_3-1}\rb$.
\end{Fact}
\begin{Fact}
(vii) Suppose the second highest $3$-valuation $\al_3-1$ occurs exactly twice.\\ Then the only possibility is to have
$\lb 3^{\al_3-1},2^2.3^{\al_3-1},2.3^{\al_3},2^23^{\al_3}\rb\subseteq X$.\\
Moreover, the solution of $X_{n=9}$ is a solution of $X_{n=7}$ ending in $\lb 2.3^{\al_3-1},3^{\al_3-1}\rb$.
\end{Fact}

\begin{Fact}
(viii) The second highest $3$-valuation $\al_3-1$ may not occur exactly thrice for a solution ending in $\lb 2.3^{\al_3},2^2.3^{\al_3}\rb$.\\
(ix) If the second highest $3$-valuation $\al_3-1$ occurs exactly thrice for a solution ending in $\lb 3^{\al_3},2.3^{\al_3}\rb$, then one of the three following situations holds:
$$\begin{array}{l}(a)\;\text{If valuation $\al_3-2$ occurs in X, then}\\\\
\al_3=4\;\text{and}\;\left\lb\begin{array}{l}\mathbf{\widehat{Z}_1}:=\lb 2,2^2,2.3^2,3^2,3^3,2.3^3,2^2.3^3,2.3^4,3^4\rb\\
\mathbf{\widehat{Z}_2}:=\lb 2,3,2.3^2,2^2.3^2,3^3,2.3^3,2^2.3^3,2.3^4,3^4\rb\end{array}\right.\;\text{are the}\\\text{only two solutions.}\\\\
\text{Moreover, $\widehat{Z}_1$ (resp $\widehat{Z}_2$) arises when the valuation $\al_3-2$ occurs exactly}\\
\text{twice in $X_{n=9}$ as $2.3^{\al_3-2}$ and $3^{\al_2}$ (resp $2.3^{\al_3-2}$ and $2^2.3^{\al_3-2}$)}\\\\
(b)\;\text{Valuation $\al_3-2$ does not occur in $X_{n=9}$ and we have a solution of $X_{n=6}$}\\
\text{ending in $\lb 2.3^{\al_3-2},2^2.3^{\al_3-2}\rb$}
\end{array}$$
\end{Fact}

\textsc{Proof of Fact $1$.} Let's treat $(i)$. Such a solution in $X_{n=9}$ would provide a solution in $X_{n=8}$ ending in $\lb 3^{\al_3-1},2^23^{\al_3-1}\rb$, which is forbidden. \\
Point $(ii)$ is straightforward. \\
Regarding $(iii)$, if such a solution existed, realize two contractions and obtain a solution in $X_{n=7}$ ending in $2.3^{\al_3-1}$ with no other occurrence with valuation $\al_3-1$. This is impossible. \\
Point $(iv)$ is similar to point $(ii)$ and is likewise straightforward. So is point $(vi)$. \\
As for $(v)$, it follows from the same argument as in $(iii)$ with $2.3^{\al_3-1}$ replaced with $3^{\al_3-1}$.\hfill $\square$\\

\textsc{Proof of Fact $2$.} If we had a solution of $X_{n=9}$ ending in this following set $\lb 3^{\al_3-1},2^2.3^{\al_3-1},2.3^{\al_3},3^{\al_3}\rb$ or that following set $\lb 3^{\al_3-1},2.3^{\al_3-1},2^2.3^{\al_3},2.3^{\al_3}\rb$ (resp $\lb 2.3^{\al_3-1},2^2.3^{\al_3-1},2.3^{\al_3},3^{\al_3}\rb$ or the set $\lb 3^{\al_3-1},2.3^{\al_3-1},2.3^{\al_3},3^{\al_3}\rb$, resp $\lb 2.3^{\al_3-1},2^2.3^{\al_3-1},2^2.3^{\al_3},2.3^{\al_3}\rb$), then we would have a solution of $X_{n=8}$ (resp $X_{n=7}$, resp $X_{n=6}$) ending in $\lb 3^{\al_3-1},2.3^{\al_3-1},2^2.3^{\al_3-1}\rb$ (resp $2.3^{\al_3-1}$, resp $3^{\al_3-1}$). None of these situations is permitted.\hfill $\square$\\

\textsc{Proof of Fact $3$.}
Suppose first the solution ends in the following set $\lb 3^{\al_3-1},2.3^{\al_3-1},2^2.3^{\al_3-1},2^2.3^{\al_3},2.3^{\al_3}\rb$.
From there, realize two contractions leading to a solution of $X_{n=7}$ containing the set $\lb 2^2.3^{\al_3-2}, 3^{\al_3-1},2^2.3^{\al_3-1}\rb$. This is impossible. \\
Suppose now the solution ends in the set $\lb 3^{\al_3-1},2.3^{\al_3-1},2^2.3^{\al_3-1},3^{\al_3},2.3^{\al_3}\rb$.
Then, by operating three contractions, we get a solution of $X_{n=6}$ which contains $\lb 2.3^{\al_3-2},2^2.3^{\al_3-2}\rb$.
The procedure is now the following. By arguments similar to those already used before, there are only three available options. Namely,
\begin{itemize}
\item[$\clubsuit$] If $\forall i\in\lb 0,1,2\rb$, we have $2^i3^{\al_3-2}\not\in X_{n=9}$, then we use induction.\\
\item[$\clubsuit\clubsuit$] If $\lb 2.3^{\al_3-2},3^{\al_3-2}\rb\subseteq X_{n=9}$, then by operating successive contractions, we obtain a solution of $X_{n=3}$ containing $2^23^{\al_3-4}$. The solutions (in non necessarily distinct integers) are well known from OEIS and are precisely $\lb 2,4,4\rb$, $\lb 2,3,6\rb$ and $\lb 3,3,3\rb$. Then the solution of $X_{n=3}$ must be $\lb 2,4,4\rb$ and the exponent $\al_3$ must be equal to $4$. Then, the original solution which we are inspecting is $\widehat{Z}_1$. \\
\item[$\clubsuit\clubsuit\clubsuit$] If $(2^2.3^{\al_3-2},2.3^{\al_3-2}\rb\subseteq X_{n=9}$, then by successive contractions, we obtain a solution of $X_{n=3}$ containing $2.3^{\al_3-3}$. Then, the solution of $X_{n=3}$ is $\lb 2,3,6\rb$ with $\al_3=3$ or $\al_3=4$ or $\lb 2,4,4\rb$ with $\al_3=3$. The latter possibility must be excluded as it leads to a solution of $X_{n=9}$ with non distinct integers. The first situation with $\al_3=3$ must also be excluded for the same reason (since then $2.3$ would appear twice). Thus, we necessarily have $\al_3=4$, leading to the solution which we denoted by $\widehat{Z}_2$.  \hfill $\square$
\end{itemize}
We will finish treating Fact $3$ $(ix)$ $(b)$ using induction and the results contained in the three facts.

First and foremost, by Fact $1$, only $2.3^{\al_3-3}$ may be the only occurrence of valuation $\al_3-3$ in $X_{n=6}$. Then, we have a solution of $X_{n=5}$ ending in $\lb 2.3^{\al_3-3},2^2.3^{\al_3-3}\rb$. By an application of the three facts, it yields a solution of $X_{n=3}$ containing either $\lb 2^2.3^{\al_3-5}\rb$ or $\lb 3^{\al_3-4},2.3^{\al_3-4}\rb$. \\
In the first case it forces $\al_3=5$ and a solution:
$$\mathbf{Z_{25}}:=\lb 2,2^2,2.3,2.3^2,3^4,2.3^4,2^2.3^4,3^5,2.3^5\rb$$
In the second situation, it also forces $\al_3=5$, hence the original solution:
$$\mbf{Z_{26}}:=\lb 2,3,2^2.3,2.3^2,3^4,2.3^4,2^2.3^4,3^5,2.3^5\rb$$

Second, valuation $\al_3-3$ could be occurring exactly twice in $X_{n=6}$, namely as $\lb 3^{\al_3-3},2^2.3^{\al_3-3}\rb$. Then, we get a solution of $X_{n=4}$ ending in $\lb 3^{\al_3-3},2.3^{\al_3-3}\rb$
Then, we have $\al_3=5$ and the original solution of $X_{n=9}$ is
$$\mbf{Z_{27}}:=\lb 2,3,3^2,2^2.3^2,3^4,2.3^4,2^2.3^4,3^5,2.3^5\rb$$

We now finish treating $(vii)$. There are several cases. \\

Assume first that $\lb 3^{\al_3-2},3^{\al_3-1},2.3^{\al_3-1}\rb\subseteq X_{n=7}$ with no other appearance of valuation $\al_3-2$ in $X_{n=9}$.
Recursively, we obtain the set of solutions:
\begin{eqnarray*}
\mbf{X_{10}}&:=&\lb 2,3,3^2,3^3,3^4,3^5,2^2.3^5,2^2.3^6,2.3^6\rb\\
\mbf{X_{11}}&:=&\lb 2,2^2,2^2.3,3^2,3^3,3^4,2^2.3^4,2.3^5,2^2.3^5\rb\\
\mbf{X_{12}}&:=&\lb 2,2^2,2.3,2^2.3^2,3^3,3^4,2^2.3^4,2.3^5,2^2.3^5\rb\\
\mbf{X_{13}}&:=&\lb 2,3,2^2.3,2^2.3^2,3^3,3^4,2^2.3^4,2.3^5,2^2.3^5\rb\\&&\\
\mbf{Y_2}&:=&\lb 2^2,3,2.3,2^2.3,3^2,3^3,2^2.3^3,2.3^4,2^2.3^4\rb
\end{eqnarray*}
Assume next that $\lb 2^2.3^{\al_3-2},3^{\al_3-1},2.3^{\al_3-1}\rb\subseteq X_{n=7}$ with no other appearance of valuation $\al_3-2$ in $X_{n=9}$. Recursively again, we obtain the other set of solutions:
\begin{eqnarray*}
\mbf{X_{14}}&:=&\lb 2,2^2,2.3,2.3^2,2^2.3^3,3^4,2^2.3^4,2.3^5,2^2.3^5\rb\\
\mbf{X_{15}}&:=&\lb 2,3,2^2.3,2.3^2,2^2.3^3,3^4,2^2.3^4,2.3^5,2^2.3^5\rb\\
\mbf{X_{16}}&:=&\lb 2,3,3^2,2^2.3^2,2^2.3^3,3^4,2^2.3^4,2.3^5,2^2.3^5\rb
\end{eqnarray*}

\noindent Assume last that $\lb 3^{\al_3-2},2.3^{\al_3-2},2^2.3^{\al_3-2},3^{\al_3-1},2.3^{\al_3-1}\rb\subseteq X_{n=7}$. Then we get a solution of $X_{n=3}$ containing $2^2.3^{\al_3-4}$. Thus $\al_3=4$ and the original solution of $X_{n=9}$ is:
$$\mbf{X_{17}}:=\lb 2,2^2,3^2,2.3^2,2^2.3^2,3^3,2^2.3^3,2.3^4,2^2.3^4\rb$$


For the processing of Fact $1$ $(ii)$, $(iv)$, $(vi)$, we will use the following automaton summarizing all the possibilities described in the three facts, together with induction.
\begin{center}
\epsfig{file=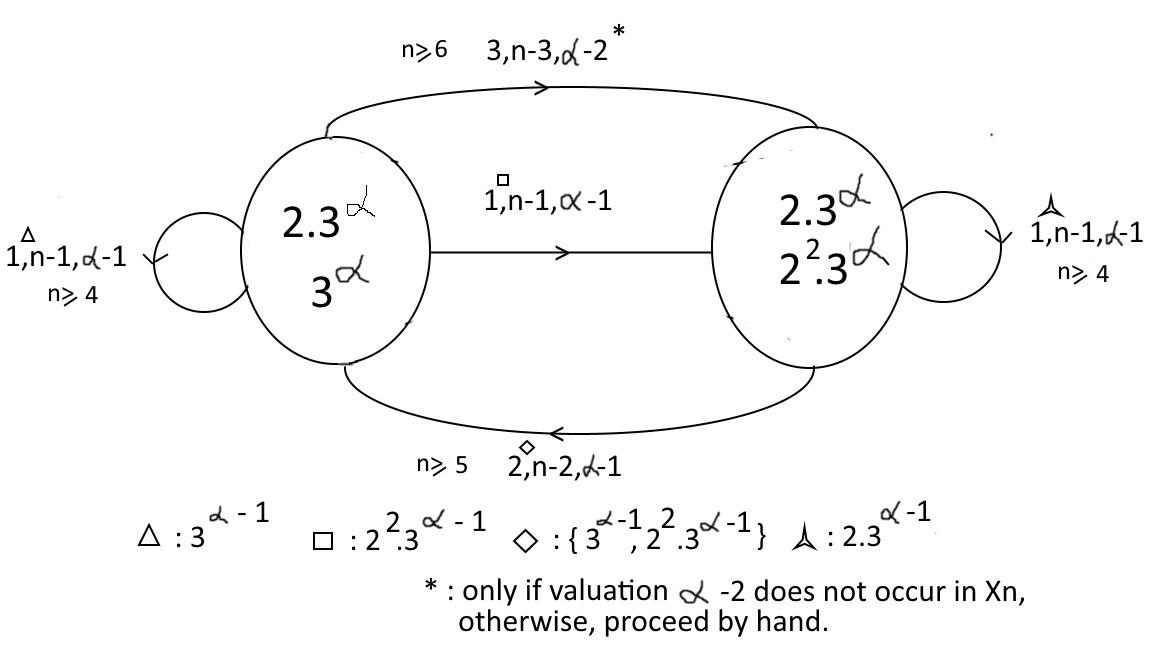, height=6cm}
\end{center}
$$\begin{array}{l}\end{array}$$
Like indicated by the star on the legend of the picture, the automaton is incomplete. There are a few leftover cases which request a special treatment. We treat them below. \\

Suppose the initial state is $\lb 2.3^{\al_3},3^{\al_3}; n=9\rb$ and we make a transition to the state $\lb 2.3^{\al_3-1},3^{\al_3-1};n=8\rb$. Suppose the solution of $X_{n=8}$ is of the form
$$\lb x_1,x_2,x_3,3^{\al_3-2},2.3^{\al_3-2},2^2.3^{\al_3-2},2.3^{\al_3-1},3^{\al_3-1}\rb$$
with valuation $\al_3-3$ occurring exactly twice. \\
Suppose first the latter valuation occurs as $3^{\al_3-3}$ and $2.3^{\al_3-3}$ in $X_{n=9}$. \\And so, the original solution is:
$$\lb x_1,3^{\al_3-3},2.3^{\al_3-3},3^{\al_3-2},2.3^{\al_3-2},2^2.3^{\al_3-2},3^{\al_3-1},3^{\al_3},2.3^{\al_3}\rb$$
By successive contractions, we get a solution of $X_{n=3}$ of the form:
$$\lb x_1,2.3^{\al_3-4},2^2.3^{\al_3-4}\rb$$
Then, $\al_3=4$ and $x_1=2^2$. It leads to the solution $\widehat{T}_1$ of Lemma $4$. \\
The other possibility is to have:
$$\lb x_1,2.3^{\al_3-3},2^2.3^{\al_3-3},3^{\al_3-2},2.3^{\al_3-2},2^2.3^{\al_3-2},3^{\al_3-1},3^{\al_3},2.3^{\al_3}\rb,$$
which by contraction leads to a solution of $X_{n=3}$ of the form:
$$\lb x_1,3^{\al_3-3},2.3^{\al_3-3}\rb$$
It forces $\al_3=4$ and $x_1=2$. We therefore obtain a new solution, namely
$$\mbf{\widehat{Z}_3}:=\lb 2,2.3,2^2.3,3^2,2.3^2,2^2.3^2,3^3,3^4,2.3^4\rb$$
Suppose now we do another transition of the same kind a second time, followed by a transition by $3$. And so, we have a tentative solution:
$$\lb x_1,x_2,3^{\al_3-3},2.3^{\al_3-3},2^2.3^{\al_3-3},3^{\al_3-2},3^{\al_3-1},3^{\al_3},2.3^{\al_3}\rb$$

\noindent After inspection by contraction, it is impossible to have
$$\lb x_1,x_2\rb=\lb 2.3^{\al_3-4},3^{\al_3-4}\rb\;\;\text{nor to have}\;\;\lb x_1,x_2\rb=\lb 2.3^{\al_3-4},2^2.3^{\al_3-4}\rb$$


\noindent These were the only additional possibilities to explore with an initial state of the form $\lb 3^{\al_3},2.3^{\al_3}\rb$. If the initial state is rather $\lb 2.3^{\al_3},2^2.3^{\al_3}\rb$, then we must transit by $2$ first and after inspection, there is no complementary addition to $X_{17}$. \\
All the additional possibilities have now been investigated, which closes the discussion.\\



From there, the remaining solutions can be read out from the two trees of appendices A and B with respective roots $\lb2,3^{\al_3},3^{\al_3}\rb$ or $\lb 2.3^{\al_3},2^2.3^{\al_3}\rb$.\\ On the respective trees, we imposed that the first transition on the automaton is by $1$, as the cases of transitions by $2$ or $3$ corresponding to $(vii)$ and $(ix)(b)$ were already processed. From the trees, it suffices to follow the path from the root to a given leaf represented by a square in order to get a solution. More explicitly, on the picture view the sets written inside the leaves as triples and conserve the orders in which the integers appear. If the father of a leaf has $4$ (resp $5$, resp $6$) vertices, then take the first two elements (resp the first element, resp the first element determined uniquely by $2^2$ in that case) of the triple, then follow the path up, finally add the root. The solution $Z_1$ is the only solution for which $\al_2=1$. This ends the proof of Lemma $4$. \hfill $\square$\\

Below the trees in respective appendices A and B, we listed the solutions in the order in which the leaves appear from left to right. The solutions $Z_{24}-Z_{26}$ (resp $X_{10}-X_{17}$) which were listed separately correspond to the leaves of the missing right hand side of the tree with root $\lb 2.3^{\al},3^{\al}\rb$ (resp $\lb 2.3^{\al},2^2.3^{\al}\rb$). Moreover, these solutions are listed and numbered respecting the same order as before, to be consistent. \\

We finally deduce Lemma $3$ from the elements of proof for Lemma $4$. \\

First and foremost, we must have $\lb 2.q^{\al_q},q^{\al_q}\rb\subseteq X$ and $q|\sigma_1(1,2)=1+2$. That is $q|3$. Then $q=3$. So the solution must end in $\lb 2.3^{\al_3},3^{\al_3}\rb$. We then use the automaton with such an initial state. But since $2^2$ is not a factor, the only transition that is allowed is the one getting to the same state again with the exponent decreased by $1$. Proceeding inductively, we get a unique solution, namely $Z_1$. \hfill $\square$\\

In order to allow for generalization when the number of vertices is $k$ with $k>9$, we summarized the moves corresponding to $\clubsuit\clubsuit$ and $\clubsuit\clubsuit\clubsuit$ in the proof of Fact $3$ into a second automaton. It reads:
\begin{center}
\epsfig{file=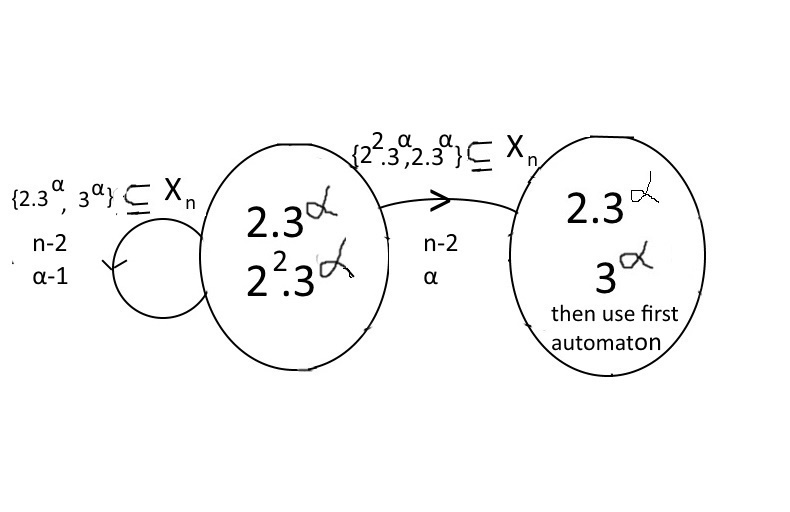, height=7cm}
\end{center}

The next part discusses the case when the number of vertices is greater than $9$, based on the general analysis done in the section above, under the same hypotheses as before and in particular under the restriction $\al_2\leq 2$.

\section{Programming with more unknowns.}

Consider $(E)_n$ in distinct integers with $n>9$. Given a solution to $(E)_n$ in $x_i$'s, suppose $2$ is a divisor of the $x_i$'s with $\al_2=2$. Suppose $q$ is an odd divisor of the $x_i$'s and so the highest power of $q$ appears twice or thrice.
\begin{itemize}
\item[$\star$] If the highest power of $q$ appears thrice, the same reasoning as in the proof of Lemma $4$ applies to a larger number of unknowns and leads to $q=7$ with the unique solution when $n$ is odd:
    $$\mathbf{V_{n\,\textbf{odd}}}:=\big\lb 2,2^2,7,2.7,7^2,2.7^2,\dots,7^{\frac{n-5}{2}},2.7^{\frac{n-5}{2}},7^{\frac{n-3}{2}},2.7^{\frac{n-3}{2}},2^2.7^{\frac{n-3}{2}}\big\rb$$
When $n$ is even, there is no solution.
\item[$\star\star$] If the highest power of $q$ appears exactly twice, then either $q=3$ with a solution ending in $\lb 2.3^{\al_3},2^2.3^{\al_3}\rb$ or in $\lb 3^{\al_3},2.3^{\al_3}\rb$, or $q=5$ with a solution necessarily ending in $\lb 5^{\al_5},2^2.5^{\al_5}\rb$. We first treat the latter case. Like before, we know that the second largest $5$-valuation occurring must be $\al_5-1$ and it must occur as $5^{\al_5-1}$. This time, working inductively, we see that there is a unique solution for each $n$, even or odd, and it reads:
    $$\mathbf{U_n}:=\lb 2,2^2,5,5^2,\dots,5^{n-4},5^{n-3},2^2.5^{n-3}\rb$$
    When $q=3$, we wrote a CAML program which returns the list of all the solutions and also counts them. This program is based on implementing the two automata of $\S\,3$. It uses a recursive transition function from a state to another state of the automa. As part of this transition function we build inductively a list that keeps track of the solution under construction every time the function gets applied recursively (these additions correspond to the labels on the edges of the trees in Appendices A and B) and the number of integer variables gets decreased accordingly.
    Ultimately, each list will be a solution to the equation. \\
    At each step after reaching a new state (corresponding to a node in the tree), the program applies all the possible transition orders. Thus, in the end we obtain all the solutions. By the way the recursion works, the tree gets traversed in depth first. \\
    Because we proceed from a state in the highest $3$-valuation whose value gets only determined after reaching a leaf in the tree, the transition function takes in argument a function which returns the
    list under construction depending on the highest $3$-valuation. The transition function also takes in argument a counter which keeps track of the variation of the exponent inside the list. \\
    We reach a leaf in the tree when the number of variables left and provided as entry in the transition function reaches a lower bound for a given order. In that case, the transition function adds the very last elements to the list output of the entry function which results into a new function, determines what the highest $3$-valuation should be depending on the value of the counter, then applies the new function with that value and returns the print-out of the completed list, that is a solution to the equation. It also increments a counter in order to keep track of the number of solutions, which also gets returned at the end of the procedure.

    Finally we briefly discuss about some of the orders appearing in the program.
    We recall that the second automaton gets used whenever a solution no longer ends in state one or two after a transition "Starone" on the first automaton. The orders "Startwo" and "Starthree" model the two possible jumps from the first automaton onto the second automaton. As far as the orders "RedStartwo" and "RedStarthree", these may occur only after a Startwo jump took place, or following a RedStartwo move.
    Namely, when you are in state two on the second automaton, you can either decide to jump back onto the first automaton at the cost of an exponent jump in the solution, or decide to stay on the second automaton, in which case there won't be any exponent jump in the solution.

    We ran the main program successfully up to $17$ integer variables, though we displayed the results in Appendix E only up to $n=13$. If we only count the solutions like in Appendices C or D, the program runs fine up to $n=30$.\\ The next section addresses the theoretical count of the solutions.
\end{itemize}

\section{Combinatorics of the number of solutions when the integers are distinct.}
When the integers are distinct, the only admissible primes are $q=3,5,7$. When $q\in\lb 5,7\rb$, the number of solutions is well determined, as in $\S\,4$. When $q=3$, we have seen that each solution is uniquely determined by a path followed from the root to a leaf on a decision tree with root either $\lb 2.3^{\al},2^2.3^{\al}\rb$ or $\lb 3^{\al},2.3^{\al}\rb$, where $\al$ is the highest $3$-valuation of the solution. Call such respective trees $T_1$ and $T_2$. Our first step in counting the solutions consists of determining the number of nodes of each tree at a given depth that is less than the minimal depth at which the first leaf in the tree appears. We denote the respective numbers by $N_i(d)$ with $i=1,2$ at a given admissible depth $d$.
The number of nodes equals the number of edges at a given depth. Below are the rules governing the number of edges. These rules are read out of the automata. The root is at depth zero and there is zero edge.

\begin{eqnarray}
0&\lra& 2\;\;(i=1),\;\; 0\;\lra\;5\;\;(i=2)\\
2&\lra& 2+5\\
5&\lra& 5+2+2+2+\h+5\\
\h&\lra& 5+2+\h
\end{eqnarray}

The rules are best illustrated on the tree of Appendix F (done when $i=1$).
We derive:
\begin{Proposition}
$N_i(d)=2a_2(d)+5a_5(d)+2a_{\h}(d)$ with:

$$\left\lb\begin{array}{ccc}
a_2(d)&=&a_2(d-1)+3a_5(d-1)+a_{\h}(d-1)\\
a_{\h}(d)&=&a_5(d-1)+a_{\h}(d-1)\\
a_5(d)&=&a_2(d-1)+2a_5(d-1)+a_{\h}(d-1)
\end{array}\right.$$

and initial conditions:
$$\begin{array}{ccc}
a_2(1)=\begin{cases} 1&\text{if $i=1$}\\
0&\text{if $i=2$}\end{cases}&
a_5(1)=\begin{cases}
0&\text{if $i=1$}\\
1&\text{if $i=2$}\end{cases}&a_{\h}(1)=0\end{array}$$
\end{Proposition}

These recurrences can be solved using Mathematica. \\

We obtain a bound for the number of solution $\ms(n)(q=3)$:
$$N_1(d^{(1)}_{min}(n))+N_2(d^{(2)}_{min}(n)\leq\ms(n)(q=3)\leq N_1(d^{(1)}_{max}(n))+N_2(d^{(2)}_{max}(n))$$

We derive the following proposition.

\begin{Proposition} Bounding the number $\ms(n)$ of solutions to $(E)_n$ in distinct integers. Let $N_1$ and $N_2$ be defined like in Proposition $2$.  \\\\
If $n\equiv 0\pmod{5}$, then $$N_1\Big(\frac{n}{5}\Big)+N_2\Big(\frac{n}{5}\Big)+\leq \ms(n)-\delta_{\text{n is even}}-2\delta_{\text{n is odd}}\leq N_1(n-3)+N_2(n-3)$$
If $n\equiv 1\pmod{5}$, then $$N_1\Big(\frac{n+4}{5}\Big)+N_2\Big(\frac{n-1}{5}\Big) \leq\ms(n)-\delta_{\text{n is even}}-2\delta_{\text{n is odd}}\leq N_1(n-3)+N_2(n-3)$$
If $n\equiv 2\pmod{5}$, then $$N_1\Big(\frac{n+3}{5}\Big)+N_2\Big(\frac{n+3}{5}\Big) \leq \ms(n)-\delta_{\text{n is even}}-2\delta_{\text{n is odd}}\leq N_1(n-3)+N_2(n-3)$$
If $n\equiv 3\pmod{5}$, then $$N_1\Big(\frac{n+2}{5}\Big)+N_2\Big(\frac{n-3}{5}\Big) \leq \ms(n)-\delta_{\text{n is even}}-2\delta_{\text{n is odd}}\leq N_1(n-3)+N_2(n-3)$$
If $n\equiv 4\pmod{5}$, then $$N_1\Big(\frac{n+6}{5}\Big)+N_2\Big(\frac{n+1}{5}\Big)\leq \ms(n)-\delta_{\text{n is even}}-2\delta_{\text{n is odd}}\leq N_1(n-3)+N_2(n-3)$$
\end{Proposition}

\textsc{Proof of Proposition $3$}. The maximal depth of reaching a leaf in the tree is given by a path in the tree of maximal length. Such a path is obtained when the number of unknowns is decreased by only one at each step. This is for instance achieved with a repeated Asterisk operation in $T_1$ and with a repeated Triangle operation in $T_2$. At depth $d$, there are $n-d$ unknowns and the process stops when $n-d=4$, that is $d=n-4$. Then $d^{(1)}_{max}(n)=d^{(2)}_{max}(n)=n-3$.  As for the minimal depth of reaching a leaf, it is provided by the shortest path to a leaf in the tree. We must distinguish between $T_1$ and $T_2$. We start with $T_1$. Our best first order is to use the Diamond move. After that we are in State $1$ and we may use the Starthree move resulting in a decrease of unknowns by $5$ as long as we no longer can.
\begin{itemize}
\item[*] If $n-2-5(d-1)=8$, apply Starthree a last time and reach a leaf. We must have $n\equiv 0\pmod{5}$ and $d^{(1)}_{min}(n)=\frac{n}{5}$.\\
\item[**] If $n-2-5(d-1)=7$, Starthree is no longer available and the best possible option is the Starone order. We must have $n\equiv 4\pmod{5}$ and we get $d_{min}^{(1)}(n)=\frac{n+6}{5}$.\\
\item[***] If $n-2-5(d-1)=6$, proceed likewise. We must have $n\equiv 3\pmod{5}$ and we obtain $d^{(1)}_{min}(n)=\frac{n+2}{5}$. \\
\item[****] If $n-2-5(d-1)=5$, the last Starthree order should be replaced with a Startwo order in order to be positioned in State $2$ and to proceed further on the path itinery with a Diamond order. We must have $n\equiv 2\pmod{5}$ and we get $d^{(1)}_{min}(n)=\frac{n+3}{5}$. \\
\item[*****] Last, if $n-2-5(d-1)=4$, we have a leaf and no further action is needed. We must have $n\equiv 1\pmod{5}$ and we get $d^{(1)}_{min}(n)=\frac{n+4}{5}$.
\end{itemize}
When processing $T_2$, the only difference with before is that we apply the Starthree moves straightaway. Then $n-2-5(d-1)$ must be replaced with $n-5d$. Respecting the same order of treatment as before, we obtain the five following possibilities:
$$\begin{array}{l}\Big(n\equiv 3\pmod{5}, d^{(2)}_{min}(n)=\frac{n-3}{5}\Big),
\Big(n\equiv 2\pmod{5}, d^{(2)}_{min}(n)=\frac{n+3}{5}\Big), \\
\Big(n\equiv 1\pmod{5}, d^{(2)}_{min}(n)=\frac{n-1}{5}\Big),
\Big(n\equiv 0\pmod{5}, d^{(2)}_{min}(n)=\frac{n}{5}\Big), \\
\Big(n\equiv 4\pmod{5}, d^{(2)}_{min}(n)=\frac{n+1}{5}\Big)\end{array}$$
Gathering all the results for both trees leads to the statement of Proposition $3$.

In fact, we can precisely count the number of solutions using a Mathematica program. It suffices to keep track of how the number of unknowns varies with the different orders Triangle, Square, Diamond, Asterisk, Starone, Startwo, Starthree, RedStartwo, RedStarthree, respectively abbreviated by tri, sq, dia, ast, St1, St2, St, t and p.

\begin{Proposition}

The following recursions hold:
\begin{eqnarray*}
ast(n)&=& ast(n-1)+ dia(n-1)+\delta_{n=4}\\
tri(n)&=& tri(n-1)+ sq(n-1)+ St1(n-1)+ St2(n-1)+ St(n-1)\\
sq(n)&=& ast(n-1)+ dia(n-1)\\
dia(n)&=& tri(n-2)+ sq(n-2)+ St1(n-2)+ St2(n-2)+ St(n-2)\\
St1(n)&=& ast(n-3)+ dia(n-3)+ \delta_{n=6}\\
St2(n)&=& ast(n-5)+ dia(n-5)+ t(n-5)+ p(n-5)+\delta_{n=8}\\
St(n)&=& tri(n-5)+ sq(n-5)+ St1(n-5)+ St2(n-5)+ St(n-5)\\
t(n)&=& ast(n-2)+ dia(n-2)+ t(n-2)+p(n-2)+ \delta_{n=5}\\
p(n)&=& tri(n-2)+ sq(n-2)+ St1(n-2)+ St2(n-2)+ St(n-2)
\end{eqnarray*}
with initial conditions:
$$\left|\begin{array}{l}dia(3)= sq(3)= ast(3)= St(3)= St2(3)= St1(3)= t(3)= p(3)= dia(4)\\
= St(4)= St2(4)= St1(4)= t(4)= p(4)= St(5)= St2(5)= St1(5)= \\
St2(6)= St(6)= St2(7)= St(7)= 0\\
tri(3)= 1
\end{array}\right.$$
The system of recurrence equations has a unique solution in $n$\\
$$(tri(n), sq(n), ast(n), dia(n), St1(n), St2(n), St(n), t(n), p(n))$$ Moreover, the sum of the first $7$ components of the $9$-tuple is precisely the number of solutions $\ms(n)(q=3)$.
\end{Proposition}
\textsc{Proof of Proposition $4$.}
The left hand sides represent the order you apply with $n$ unknowns. The right hand side gathers all the orders you may operate after reduction next, with the adequate decreased number of unknowns. If there are not enough unknowns left, it will be impossible to apply a given order and so its contribution to the number of solutions is null. This is conveyed through the initial conditions with equalities to zero. We must now specify the other initial conditions.
\begin{itemize}
\item[1.] When $n=3$, Triangle is the only admissible order and it leads to a solution. Hence $tri(3)=1$.
\item[2.] When $n=4$, Triangle and Asterisk are the only admissible orders, each leading to a solution. We see that we must add a corrective term in the first recurrence equation in order to set $ast(4)=1$.
\item[3.] When $n=5$, RedStartwo, RedStarthree, Asterisk, Triangle, Square, Diamond are all admissible orders, amongst which only RedStartwo needs a correction in order to have $t(5)=1$.
\item[4.] When $n=6$, we must add Starone to the list of admissible orders and it needs a correction in order to have $St1(6)=1$.
\item[5.] When $n=7$, no new order gets added.
\item[6.] When $n=8$, both Startwo and Starthree get finally added. Only Startwo does need a correction in order to get $St2(8)=1$.
\end{itemize}

The recurrence system has the same order as the number of initial conditions, hence has a unique solution.
In order to find out the number of solutions to $(E)_n$ when $q=3$, we start with $n$ variables and apply all of Triangle, Square, Starone, Startwo, Starthree, Asterisk, Diamond. This ends the proof.\hfill $\square$

\section{Concluding words and prospects.}

It would be interesting to generalize the results by offering an algorithm for finding all the distinct solutions of $(E)_n$ with only two prime divisors and of the form $2^{\al}q^{\beta}$. That is we now allow the highest $2$-valuation $\al_2$ to satisfy to $\al_2\geq 3$. \\

\noindent \textbf{Problem $1$. Can we count all the solutions in distinct integers of the form $2^{\al}q^{\beta}$ with $q$ any odd prime ?}\\

\noindent It is to expect that dealing with the cases $\al_2\geq 3$ will be a generalization of the work contained in the present paper, with even more combinatorics involved when it comes to counting all the solutions. \\

\noindent \textbf{Problem $2$. More generally, can we count all the solutions in distinct integers with only two prime divisors ?}\\

\noindent In the case of nine integer variables, it is known since \cite{BUR} that there are only five solutions in distinct odd integers. These were recalled at the beginning of $\S\,3$.
After scrutiny, none of $\mbf{B_1-B_5}$ is composed of only two prime divisors and the minimum number of prime divisors is $3$. Thus, in the case of $(E)_9$, if we can finish solving Problem $1$, we would get the total number of solutions in distinct integers with only two prime divisors. \\

\noindent Problem $2$ is tackled in \cite{GAO} for two or more primes in the case when the primes are set (contrary to in our approach). The authors denote by $T_n(p_1,\dots,p_t)$ these numbers and one of their main results is the existence of two constants $n_0$ and $c$ that depend on the chosen primes, such that for any $n>n_0$, we have either $T_n(p_1,\dots,p_t)=0$ or $T_n(p_1,\dots,p_t)>c^n$. They leave as an open problem whether there would exist such an upper bound. While providing a lower bound in the special case of $T_k(2,3,5)$, they state some interesting decompositions of $1/a$ as a sum of two unit fractions when $2$ divides $a$, involving the prime $3$ and four unit fractions when $2^2$ divides $a$, involving the primes $2$, $3$ and $5$. These identities namely read:
\begin{eqnarray}
\frac{1}{a}&=&\frac{1}{\frac{5a}{4}}+\frac{1}{10a}+\frac{1}{15a}+\frac{1}{30a}\\
\frac{1}{a}&=&\frac{1}{\frac{3a}{2}}+\frac{1}{3a}
\end{eqnarray}
In the case of $(E)_9$, we show how we can use these identities to derive solutions in $2^{\al}3^{\beta}5^{\gamma}$.
We expand a solution of $(E)_3$ instead of retracting a solution to one of $(E)_3$ like we did in $\S\,3$.
First and foremost, we have:
\begin{equation}
\frac{1}{2}+\frac{1}{2^2}+\frac{1}{5}+\frac{1}{5.2^3}+\frac{1}{3.5.2^2}+\frac{1}{3.5.2^3}=1
\end{equation}
by starting from the solution $\lb 2,4,4\rb$ of $(E)_3$ and using $(17)$ with $a=4$. We obtain a solution of $(E)_6$. From there, each unit fraction containing $2^2$ (resp $2$) as a divisor may be expanded using Identity $(17)$ (resp Identity $(18)$), which yields four solutions in $2^{\al}3^{\beta}5^{\gamma}$ (resp $5$ solutions of $(E)_7$ and at least $5^3$ solutions of $(E)_9$ since $3a$ remains divisible by $2$ when $a$ is divisible by $2$).\\
The authors of \cite{GAO} also study the special case $T_n(3,5,7)$. They show that:
$$T_n(3,5,7)\geq C^{n}\sqrt{62}\qquad\forall n\geq 11$$
In \cite{BUR2}, Burshtein showed that $T_{11}(3,5,7)=17$. By confronting both works, we deduce:
$$C\leq e^{\frac{1}{11}(ln(17)-\frac{1}{2}ln(62)}\simeq 1.072$$

\noindent Finally, we note that the questions we raise and those raised in \cite{GAO} are both different from the almost century old questions raised by Erd\"os. He and Graham had been investigating Egyptian fractions with each denominator having three distinct prime divisors but while doing so they allow the three primes to vary from a unit fraction to the other, cf. eg. \cite{GER}. Their work led to a paper \cite{BEG} involving also Butler in $2015$. Erd\"os is a co-author though he died twenty years earlier in $1996$. They show that each natural number has a representation in Egyptian fraction with each denominator having three distinct prime divisors. In \cite{BUR3}, Burshtein also discusses the case of two distinct prime divisors. \\

\textsc{\textbf{Acknowledgements.}} Thanks go to Joel Louwsma whose interesting work on arithmetical structures originated the present work. His insight on the Diophantine equation in non necessarily distinct integers and completely different approach allowed the author to identify a small mistake in her original program. Joel Louwsma namely allowed to confront, then confirm by different means than those presented here the counting results of the present paper. He is gratefully acknowledged for that. An oral presentation of the author at the Southern California Discrete Mathematics Symposium SoCalDM $2025$ at the University of California at Irvine allowed to notice during a computer demo that two solutions differed only by one integer at length $12$. This was simply due to a typo in the program. The author is indebted to Nathan Kaplan and other attendees of the conference for their scrutiny which allowed to correct the program accordingly. 
My Profound, Happy and Sincere thanks go to MF and PA for allowing to identify a missing leaf in the tree of Appendix B. Their generous time and excited dedication allowed to compare the hand solutions written in terms of powers with the computer solutions for $n=9$ written in terms of integers, thus fixing a count difference of a simple unit and allowing again some verifications on the results. I am indebted for their patience, generosity and success. Finally, my most sincere thanks go to KW for her great encouragements and very kind support. 
\newpage \begin{Large}Appendix A\end{Large}\\
$\begin{array}{cc}\nts\nts\nts\nts\nts\nts\nts\nts\nts\nts\nts\nts\nts\nts\nts\nts\nts\nts\nts\nts\nts\nts\nts\nts\nts\nts\nts\nts\nts\nts\nts\nts
\text{\epsfig{file=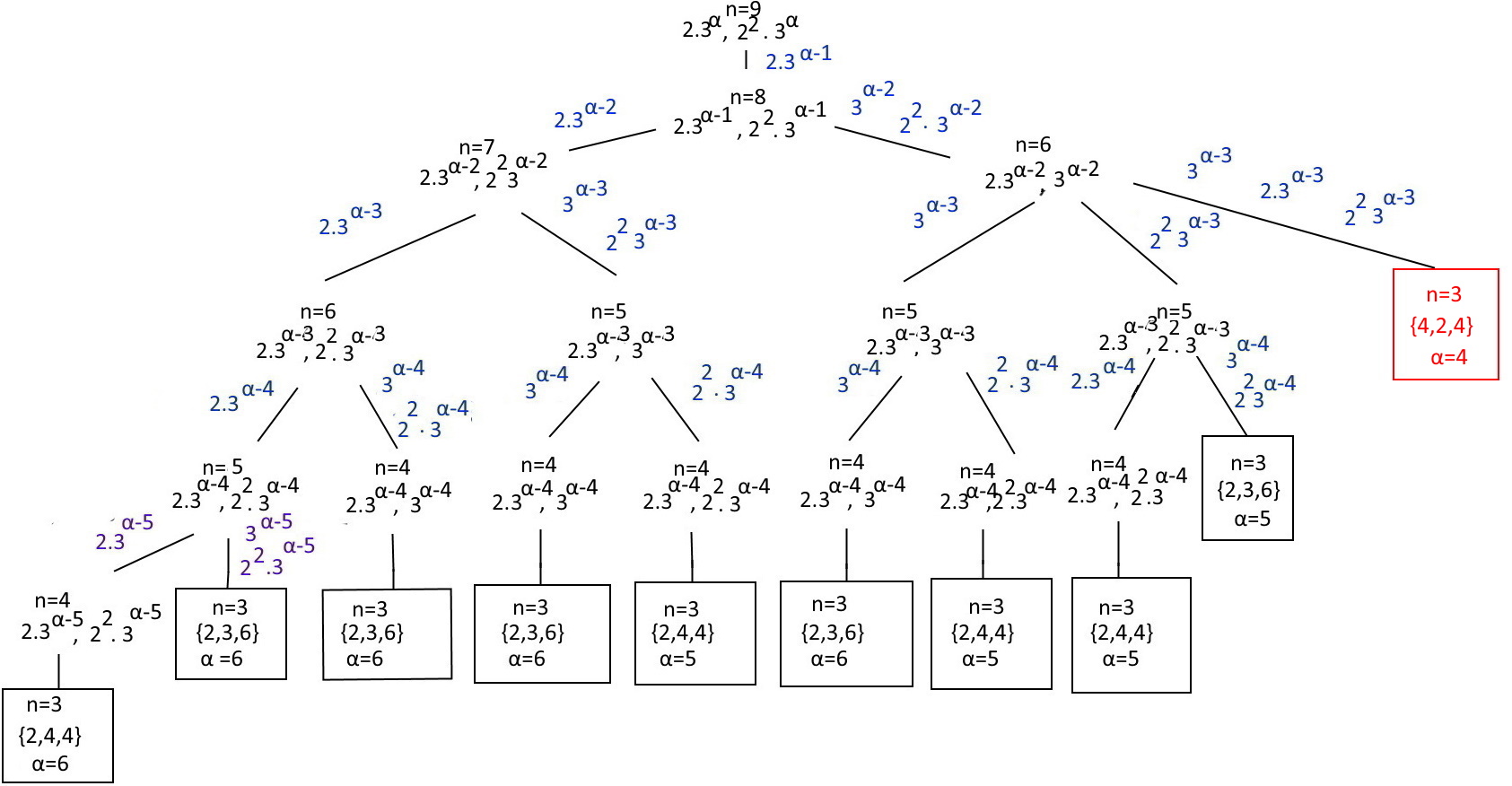, height=10cm}}&\qquad\qquad\qquad\qquad\qquad\end{array}$
\begin{eqnarray*}
\mbf{X_1}&:=&\lb 2,2^2,2.3,2.3^2,2.3^3,2.3^4,2.3^5,2.3^6,2^2.3^6\rb\\
\mbf{X_2}&:=&\lb 2,3,2^2.3,2.3^2,2.3^3,2.3^4,2.3^5,2.3^6,2^2.3^6\rb\\
\mbf{X_3}&:=&\lb 2,3,2^2.3^2,3^2,2.3^3,2.3^4,2.3^5,2.3^6,2^2.3^6\rb\\
\mbf{X_4}&:=&\lb 2,3,3^2,3^3,2^2.3^3,2.3^4,2.3^5,2.3^6,2^2.3^6\rb\\
\mbf{X_5}&:=&\lb 2,2^2,2^2.3,3^2,2^2.3^2,2.3^3,2.3^4,2.3^5,2^2.3^5\rb\\
\mbf{X_6}&:=&\lb 2,3,3^2,3^3,3^4,2^2.3^4,2.3^5,2.3^6,2^2.3^6\rb\\
\mbf{X_7}&:=&\lb 2,2^2,2^2.3,3^2,3^3,2^2.3^3,2.3^4,2.3^5,2^2.3^5\rb\\
\mbf{X_8}&:=&\lb 2,2^2,2.3,2^2.3^2,3^3,2^2.3^3,2.3^4,2.3^5,2^2.3^5\rb\\
\mbf{X_9}&:=&\lb 2,3,2^2.3,2^2.3^2,3^3,2^2.3^3,2.3^4,2.3^5,2^2.3^5\rb\\
&&\\
\mbf{Y_1}&:=&\lb 2^2,3,2.3,2^2.3,3^2,2^2.3^2,2.3^3,2.3^4,2^2.3^4\rb
\end{eqnarray*}
\newpage \begin{Large}Appendix B\end{Large}\\
$\begin{array}{cc}\nts\nts\nts\nts\nts\nts\nts\nts\nts\nts\nts\nts\nts\nts\nts\nts\nts\nts\nts\nts\nts\nts\nts\nts\nts\nts\nts\nts\nts\nts\nts\nts
\nts\nts\nts\nts\nts\nts\nts\nts\nts\nts\nts\nts\nts\nts\nts\nts\nts\text{\epsfig{file=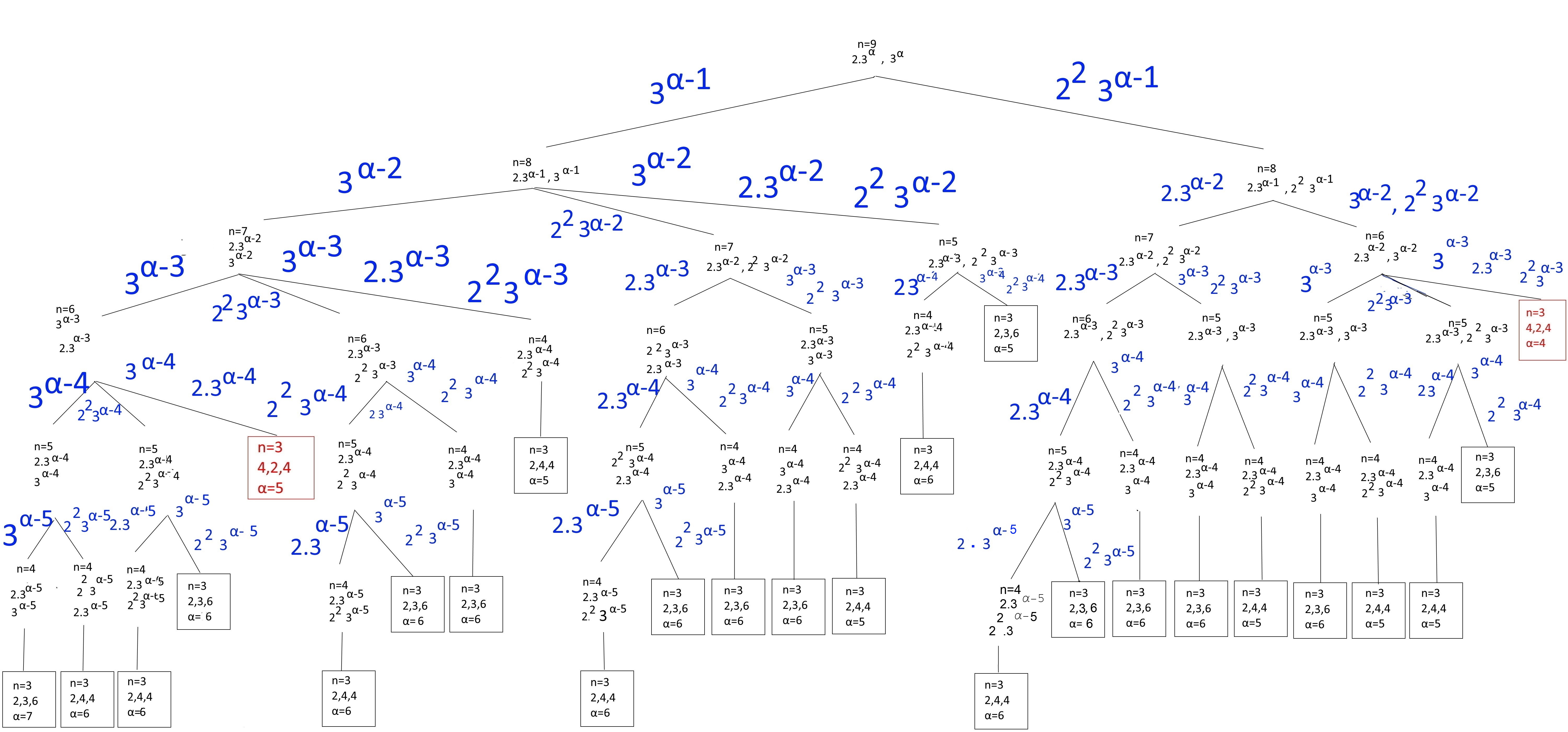, height=10cm}}&\qquad\qquad\qquad\qquad\qquad\end{array}$

\begin{eqnarray*}
\mbf{Z_1}&:=&\lb 2,3,3^2,3^3,3^4,3^5,3^6,3^7,2.3^7\rb\\
\mbf{Z_2}&:=&\lb 2,2^2,2^2.3,3^2,3^3,3^4,3^5,3^6,2.3^6\rb\\
\mbf{Z_3}&:=&\lb 2,2^2,2.3,2^2.3^2,3^3,3^4,3^5,3^6,2.3^6\rb\\
\mbf{Z_4}&:=&\lb 2,3,2^2.3,2^2.3^2,3^3,3^4,3^5,3^6,2.3^6\rb\\
\mbf{Z_5}&:=&\lb 2,2^2,2.3,2.3^2,2^2.3^3,3^4,3^5,3^6,2.3^6\rb\\
\mbf{Z_6}&:=&\lb 2,3,2^2.3,2.3^2,2^2.3^3,3^4,3^5,3^6,2.3^6\rb\\
\mbf{Z_7}&:=&\lb 2,3,3^2,2^2.3^2,2^2.3^3,3^4,3^5,3^6,2.3^6\rb\\
\mbf{Z_8}&:=&\lb 2,2^2,3^2,2.3^2,2^2.3^2,3^3,3^4,3^5,2.3^5\rb\\
\mbf{Z_9}&:=&\lb 2,2^2,2.3,2.3^2,2.3^3,2^2.3^4,3^5,3^6,2.3^6\rb\\
\mbf{Z_{10}}&:=&\lb 2,3,2^2.3,2.3^2,2.3^3,2^2.3^4,3^5,3^6,2.3^6\rb\\
\mbf{Z_{11}}&:=&\lb 2,3,3^2,2^2.3^2,2.3^3,2^2.3^4,3^5,3^6,2.3^6\rb\\
\mbf{Z_{12}}&:=&\lb 2,3,3^2,3^3,2^2.3^3,2^2.3^4,3^5,3^6,2.3^6\rb\\
\mbf{Z_{13}}&:=&\lb 2,2^2,2^2.3,3^2,2^2.3^2,2^2.3^3,3^4,3^5,2.3^5\rb\\
\mbf{Z_{14}}&:=&\lb 2,2^2,2.3,3^3,2.3^3,2^2.3^3,3^4,3^5,2.3^5\rb\\
\mbf{Z_{15}}&:=&\lb 2,3,2^2.3,3^3,2.3^3,2^2.3^3,3^4,3^5,2.3^5\rb
\end{eqnarray*}
\begin{eqnarray*}
\mbf{Z_{16}}&:=&\lb 2,2^2,2.3,2.3^2,2.3^3,2.3^4,2^2.3^5,3^6,2.3^6\rb\\
\mbf{Z_{17}}&:=&\lb 2,3,2^2.3,2.3^2,2.3^3,2.3^4,2^2.3^5,3^6,2.3^6\rb\\
\mbf{Z_{18}}&:=&\lb 2,3,3^2,2^2.3^2,2.3^3,2.3^4,2^2.3^5,3^6,2.3^6\rb\\
\mbf{Z_{19}}&:=&\lb 2,3,3^2,3^3,2^2.3^3,2.3^4,2^2.3^5,3^6,2.3^6\rb\\
\mbf{Z_{20}}&:=&\lb 2,2^2,2^2.3,3^2,2^2.3^2,2.3^3,2^2.3^4,3^5,2.3^5\rb\\
\mbf{Z_{21}}&:=&\lb 2,3,3^2,3^3,3^4,2^2.3^4,2^2.3^5,3^6,2.3^6\rb\\
\mbf{Z_{22}}&:=&\lb 2,2^2,2^2.3,3^2,3^3,2^2.3^3,2^2.3^4,3^5,2.3^5\rb\\
\mbf{Z_{23}}&:=&\lb 2,2^2,2.3,2^2.3^2,3^3,2^2.3^3,2^2.3^4,3^5,2.3^5\rb\\
\mbf{Z_{24}}&:=&\lb 2,3,2^2.3,2^2.3^2,3^3,2^2.3^3,2^2.3^4,3^5,2.3^5\rb\\&&\\
\mbf{T_1}&:=&\lb 2^2,3,2.3,2^2.3,3^2,2^2.3^2,2^2.3^3,3^4,2.3^4\rb\\
\mbf{T_2}&:=&\lb 2^2,3,2.3,2^2.3,3^2,3^3,3^4,3^5,2.3^5\rb
\end{eqnarray*}

\newpage

\noindent\begin{Large}Appendix C. Program in CAML for counting the solutions \end{Large}
\begin{verbatim}

let p n = power 3.0 n;;                            (* defines 3^n *)


type Order = Triangle | Square | Asterisk | Diamond | Starone | Startwo | Starthree
             | RedStartwo | RedStarthree ;;

                                                   (* defines the type Order *)

let rec transition (state,n,order) =

      match (state,n,order) with

            (1,n,_) when n=4 -> incre()

            |(2,n,_) when n=4 -> incre()

            |(1,n,Triangle) when n=5 ->  transition (1,4,Triangle)

            |(1,n,Square) when n=5 -> transition (2,4,Asterisk)

            |(2,n,Asterisk) when n=5  ->  transition (2,4,Asterisk)

            |(2,n,Diamond) when n=5 ->  incre()

            |(2,n,RedStartwo) when n=5 -> incre()

            |(2,n,RedStarthree) when n=5 -> incre()

            |(2,n,Diamond) when n=6 -> transition (1,4,Triangle)

            |(1,n,Starone) when n=6 -> incre()

            |(2,n,RedStartwo) when n=6 -> incre()

            |(2,n,RedStarthree) when n=6 -> transition (1,n-2,Triangle)

            |(1,n,Triangle) when n=6 -> transition (1,5,Triangle);
                                        transition (1,5,Square)

            |(2,n,Diamond) when n=7 -> transition (1,5,Triangle);
                                       transition (1,5,Square)

            |(1,n,Starone) when n=7 -> incre()

            |(2,n,RedStarthree) when n=7 -> transition (1,n-2,Triangle);
                                            transition (1,n-2,Square)

            |(1,n,Triangle) when (n=7 or n=8) -> transition (1,n-1,Triangle);
                                                 transition (1,n-1,Square);
                                                 transition (1,n-1,Starone)

            |(1,n,Startwo) when n=8 -> incre()

            |(1,n,Starthree) when n=8 -> incre()

            |(2,n,Diamond) when (n=8 or n=9) -> transition (1,n-2,Triangle);
                                                transition (1,n-2,Square);
                                                transition (1,n-2,Starone)

            |(2,n,RedStarthree) when (n=8 or n=9) -> transition (1,n-2,Triangle);
                                                     transition (1,n-2,Square);
                                                     transition (1,n-2,Starone)

            |(1,n,Startwo) when n=9 -> incre()

            |(1,n,Starthree) when n=9 -> incre()

            |(1,n,Starthree) when n=10 -> transition (1,n-5,Triangle);
                                          transition (1,n-5,Square)

            |(1,n,Starthree) when (n=11 or n=12) -> transition (1,n-5,Triangle);
                                                    transition (1,n-5,Square);
                                                    transition (1,n-5,Starone)

            |(1,n,Triangle) when n>7 -> transition (1,n-1,Triangle);
                                        transition (1,n-1,Square);
                                        transition (1,n-1,Starone);
                                        transition (1,n-1,Startwo);
                                        transition (1,n-1,Starthree)

            |(2,n,Asterisk) when n>5 -> transition (2,n-1,Asterisk);
                                        transition (2,n-1,Diamond)

            |(2,n,Diamond) when n>9 -> transition (1,n-2,Triangle);
                                       transition (1,n-2,Square);
                                       transition (1,n-2,Starone);
                                       transition (1,n-2,Startwo);
                                       transition (1,n-2,Starthree)

            |(1,n,Square) when n>5 -> transition(2,n-1,Asterisk);
                                      transition(2,n-1,Diamond)

            |(1,n,Starone) when n>7 -> transition (2,n-3,Asterisk);
                                       transition (2,n-3,Diamond)

            |(1,n,Startwo) when n>9 -> transition (2,n-5,Asterisk);
                                       transition (2,n-5,Diamond);
                                       transition (2,n-5,RedStartwo);
                                       transition (2,n-5,RedStarthree)

            |(1,n,Starthree) when n>12 -> transition (1,n-5,Triangle);
                                          transition (1,n-5,Square);
                                          transition (1,n-5,Starone);
                                          transition (1,n-5,Startwo);
                                          transition (1,n-5,Starthree)

            |(2,n,RedStartwo) when n>6 -> transition (2,n-2,Asterisk);
                                          transition (2,n-2,Diamond);
                                          transition (2,n-2,RedStartwo);
                                          transition (2,n-2,RedStarthree)


            |(2,n,RedStarthree) when n>9 -> transition (1,n-2,Triangle);
                                            transition (1,n-2,Square);
                                            transition (1,n-2,Starone);
                                            transition (1,n-2,Startwo);
                                            transition (1,n-2,Starthree)

            |(_,_,_) -> failwith "error";;

let N m = (* Returns the number of solutions in 2^a.3^b with a<3 for m variables *)
 let ct = ref 0 and cc = ref 0 and ceu = ref 0 and ced = ref 0 and cet = ref 0 and
     ca = ref 0 and cd = ref 0
      in
       let S m order =  (* with m the number of variables *)
        let c = ref 0
         in
        let incre()= c:= !(c) +1
       in
       match order with
        Triangle -> begin transition (1,m,Triangle);
                   ct:= !c end
       |Square -> begin transition (1,m,Square);
                 cc:= !c end
       |Starone -> begin transition (1,m,Starone);
                  ceu:= !c end
       |Startwo -> begin transition (1,m,Startwo);
                  ced:= !c end
       |Starthree -> begin transition (1,m,Starthree);
                    cet:= !c end
       |Asterisk -> begin transition (2,m,Asterisk);
                   ca:= !c end
       |Diamond -> begin transition (2,m,Diamond);
                  cd:= !c end
       |_ -> failwith("error")

        in

      S m Triangle; S m Square; S m Starone; S m Startwo; S m Starthree;
      S m Asterisk; S m Diamond;
      let total = !(ct) + !(cc)+ !(ceu)+ !(ced)+ !(cet)+ !(ca)+ !(cd)
         in
         print_string("There are ");print_int(total);print_string(" solutions");;

Compiling:

p : float -> float = <fun>
#Type Order defined.
#N : int -> unit = <fun>

Running the program:

#N(9);;
There are 52 solutions- : unit = ()
#N(22);;
There are 228102 solutions- : unit = ()
#N(30);;
There are 39590576 solutions- : unit = ()
#
\end{verbatim}
\newpage
\noindent\begin{Large}Appendix D. Program in CAML for counting the solutions up to a given n\end{Large}
\begin{verbatim}

let list_number_of_solutions_up_to k = (* solutions in 2^a3^b with a<3 and
                                          k the number of variables *)
  let aux m =
 let ct = ref 0 and cc = ref 0 and ceu = ref 0 and ced = ref 0 and cet = ref 0 and
     ca = ref 0 and cd = ref 0
 in
  let S m order =
   let c = ref 0
   in
  let incre()= c:= !(c) +1
   in
    match order with
       Triangle -> begin transition (1,m,Triangle);
                   ct:= !c end
      |Square -> begin transition (1,m,Square);
                 cc:= !c end
      |Starone -> begin transition (1,m,Starone);
                  ceu:= !c end
      |Startwo -> begin transition (1,m,Startwo);
                  ced:= !c end
      |Starthree -> begin transition (1,m,Starthree);
                    cet:= !c end
      |Asterisk -> begin transition (2,m,Asterisk);
                   ca:= !c end
      |Diamond -> begin transition (2,m,Diamond);
                  cd:= !c end
      |_ -> failwith("error")

        in

      S m Triangle; S m Square; S m Starone; S m Startwo; S m Starthree;
      S m Asterisk; S m Diamond;
      !(ct) + !(cc)+ !(ceu)+ !(ced)+ !(cet)+ !(ca)+ !(cd)

             in
                for n = 9 to k
                     do
                     print_string("(");print_int(n);
                     print_string(",");print_int(aux(n));print_string(")");
                     print_newline()
                     done;;

Compiling:

p : float -> float = <fun>
#Type Order defined.
#list_number_of_solutions_up_to : int -> unit = <fun>

Running the program:

#list_number_of_solutions_up_to(35);;
#(9,52)
(10,100)
(11,190)
(12,362)
(13,690)
(14,1314)
(15,2504)
(16,4770)
(17,9088)
(18,17314)
(19,32986)
(20,62844)
(21,119728)
(22,228102)
(23,434572)
(24,827932)
(25,1577348)
(26,3005110)
(27,5725234)
(28,10907522)
(29,20780642)
(30,39590576)
(31,75426626)
(32,143700256)
(33,273772866)
(34,521582802)
(35,993701908)
- : unit = ()
#
\end{verbatim}

\newpage
\noindent\begin{Large}Appendix E. Program in CAML applied to $K_n$ with
$n=9,10,11,12,13$ \end{Large}\\


\newpage
\begin{Large} Appendix F.\end{Large}
$$\begin{array}{l}\\\end{array}$$
$\begin{array}{cc}\nts\nts\nts\nts\nts\nts\nts\nts\nts\nts\nts\nts\nts\nts
\text{\epsfig{file=first_recurrence_system.jpg, height=12cm}}&\qquad\qquad\end{array}$

\newpage

\begin{Large}Appendix G. Mathematica program for counting the solutions in distinct integers.\end{Large}\\
\begin{verbatim}
In= h[k_] := RecurrenceTable[{ast[n] ==
    ast[n - 1] + dia[n - 1] + KroneckerDelta[n, 4],
   tri[n] ==
    tri[n - 1] + sq[n - 1] + St1[n - 1] + St2[n - 1] + St[n - 1],
   sq[n] == ast[n - 1] + dia[n - 1],
   dia[n] ==
    tri[n - 2] + sq[n - 2] + St1[n - 2] + St2[n - 2] + St[n - 2],
   St1[n] == ast[n - 3] + dia[n - 3] + KroneckerDelta[n, 6],
   St2[n] ==
    ast[n - 5] + dia[n - 5] + t[n - 5] + p[n - 5] +
     KroneckerDelta[n, 8],
   St[n] ==
    tri[n - 5] + sq[n - 5] + St1[n - 5] + St2[n - 5] + St[n - 5],
    t[n] ==
    ast[n - 2] + dia[n - 2] + t[n - 2] + p[n - 2] +
     KroneckerDelta[n, 5],
   p[n] ==
    tri[n - 2] + sq[n - 2] + St1[n - 2] + St2[n - 2] + St[n - 2],
   dia[4] == dia[3] == sq[3] == St[3] == St2[3] == St1[3] == St[4] ==
    St2[4] == St1[4] == St[5] == St2[5] == St1[5] == St2[6] ==
    St[6] == St2[7] == St[7] == t[4] == p[4] == t[3] == p[3] ==
    ast[3] == 0, tri[3] == 1}, {tri, sq, ast, dia, St1, St2, St, t,
   p}, {n, k, k}]; h[35]

Out= {{330140577, 191442225, 191442225, 173287025, 52743872, 29586769,
  25059215, 204595521, 173287025}}

In= 330140577 + 191442225 + 191442225 + 173287025 + 52743872 + 29586769 +
25059215

Out= 993701908
\end{verbatim}

\end{document}